\documentclass[12pt]{amsart}

\usepackage{amsmath,amssymb}

\newtheorem{thm}{Theorem}
\newtheorem{prop}{Proposition}
\newtheorem{cor}{Corollary}
\newtheorem{lemma}{Lemma}

\theoremstyle{definition}
\newtheorem{defn}{Definition}

\newtheorem{rem}{Remark}

\newcommand{\C}{{\mathbb{C}}}
\newcommand{\Z}{{\mathbb{Z}}}
\newcommand{\F}{{\mathbb{F}}}
\newcommand{\zero}{{\mathbf{0}}}

\DeclareMathOperator{\tr}{tr}
\DeclareMathOperator{\Tr}{Tr}

\begin{document}

\title[Elliptic curves and irreducible polynomials]{Elliptic curves and
explicit enumeration of irreducible polynomials with two coefficients
prescribed}

\author{Marko Moisio}
\email{mamo@uwasa.fi}
\address{Department of Mathematics and Statistics, University of Vaasa, P.O.
Box 700, FIN-65101 Vaasa, Finland}
\author{Kalle Ranto}
\email{kara@utu.fi}
\address{Department of Mathematics, University of Turku, FIN-20014
Turku, Finland}
\thanks{The work of the second author is funded by the
Academy of Finland, grant 108238}

\begin{abstract}Let $\F_q$ be a finite field of characteristic $p=2,3$. We
give the number of irreducible polynomials $x^m+a_{m-1}x^{m-1}+\cdots+a_0\in\F_q[x]$ with $a_{m-1}$ and $a_{m-3}$ prescribed for any given $m$ if $p=2$, and with $a_{m-1}$ and $a_1$ prescribed for $m=1,\dots,10$ if $p=2,3$.  
\end{abstract}


\maketitle

\section{Introduction}

Let $p$ be a prime, let $r,m,k$ a positive integers, and let $\F_q$ denote the
finite field with $q=p^r$ elements. The determination of the number $N$ of
irreducible polynomials 
$$
x^m+a_{m-1}x^m+\cdots+a_1x+a_0\in\F_q[x],
$$
with $k$ of the coefficients $a_0,\dots,a_{m-1}$ prescribed, is a difficult
problem in general and has been a subject of study for a long time, see
e.g.~\cite[p.~340]{Cohen05}, where a short survey on recent results on this
topic is given. For example, in the case of fixed $a_{m-1}$ and $a_0$ the
number $N$ was obtained by Carlitz in~\cite{Carlitz52} and later,
elementarily, by Yucas in~\cite{Yucas06}. 

In this paper the following special cases are considered:

\begin{enumerate}
\item[(i)] $p=2$, $a_{m-1}=0$, $a_{m-3}=c$,
\item[(ii)] $p=2$ or $p=3$, $a_{m-1}=c$, $a_1=0$, 
\end{enumerate}

with $c$ a given element from $\F_q$.

In these cases the problem of determining $N$ is closely related to the
problem of counting the number of rational points on the fibre products of
certain super-singular elliptic curves and of certain Kloosterman curves over
$\F_{q^m}$. This problem can be tackled by using some properties of
Kloosterman sums and cubic exponential sums, properties of Dickson
polynomials, and the well-known (see e.g.~\cite{Schoof95,Schoof91,Geer92})
weight distributions of the dual of the binary two-error-correcting BCH code
of length $q-1$ and of the binary and ternary Melas codes of length $q-1$.
This approach enables us to give $N$ explicitely for any $m$ in case (i) and
for $m=1,\dots,10$ in case (ii). 

The rest of this note is organized as follows. In Section 2 we first recall
some basic properties of Dickson polynomials and then a formula expressing the
number of rational points on certain Artin--Schreier curves in terms of
exponential sums is derived. That formula is used in Section 3 to get
explicitely the number of rational points on the fibre product of certain
super-singular elliptic curves, and finally, in Section 4, the number $N$ of
irreducible polynomials is determined by connecting it to the number of
rational points on the curves studied in Sections 2 and 3. 

\section{Preliminaries}

In this section some notations are fixed, a result concerning the point
counting on fibre products of certain Artin--Schreier curves is established,
and some results from~\cite{Lidl93, Moisio06b,Moisio06c} are recalled.

Let $\Tr$ an $\tr$ denote the trace functions from $\F_{q^m}$ onto $\F_p$ and
$\F_q$, respectively. 

Let $\omega$ be a complex number, and let 
$$
D_m(T,\omega):=\sum_{j=0}^{\lfloor
m/2\rfloor}\frac{m}{m-j}\binom{m-j}{j}(-\omega)^j
T^{m-2j}\in\C[T].
$$
denote the Dickson polynomial of the first kind of degree $m$ with parameter
$\omega$. We shall need the following two fundamental properties of Dickson
polynomials: 

\begin{equation}\label{e:Dvalues}
D_m(t,\omega)=\left(\frac{t+\sqrt{t^2-4\omega}}2\right)^m+
         \left(\frac{t-\sqrt{t^2-4\omega}}2\right)^m\quad\forall\,t\in\C,
\end{equation}
\begin{equation}\label{e:Dfunct}
D_m\Bigl(t+\frac{\omega}t,\omega\Bigr)=t^m+\frac{\omega^m}{t^m}\quad\forall\,
t\in\C^*.
\end{equation}

Next we consider point counting on fibre products of Artin--Schreier
curves. Let $L$ be a $\F_{p}$-subspace of the rational function field
$\F_{q^m}(x)$ with a basis $\{f_1,\dots,f_n\}\subset\F_{q^m}(x)\setminus\F_q$, and
assume that the multiplicities of poles of the non-zero elements of $L$ are
not divisible by $p$.

Let $f(x)\in L\setminus\{\zero\}$, let $F_{f,m}=\F_{q^m}(x,y_f)$ with
$y_f^p-y_f=f(x)$, and let 
$$
F_m=\F_{q^m}(x,y_{f_1},\dots,y_{f_n}).
$$  
Let $e$ be the canonical additive character of $\F_{q^m}$.

\begin{prop}[{\cite[Theorem 8]{Moisio06c}}]\label{t:expsumplaces}The
number $N_m$ of rational places of $F_m$ is given by 
$$
N_{m}=q^m+1-\frac{p^n-1}{p-1}+\sum_{f\in L\setminus\{0\}/\F_p^*}|S_f|
+\sum_{f\in L\setminus\{0\}}\sum_{z\in\F_{q^m}\setminus P_f}e(f(z)),
$$
where $S_f$ is the set of rational places of $F_{f,m}$ lying above
$P_{\infty}$, and $P_f$ is the set of poles of $f(x)$ in $\F_{q^m}$.
\end{prop}

In~\cite{Moisio06c} Proposition~\ref{t:expsumplaces} was applied to the
fibre product $L_m$ of Kloosterman curves defined by
$$
L_m=\F_{q^m}(x,y_a,y_b),\ y_a^q-y_a=x+ax^{-1},\ y_b^q-y_b=x+bx^{-1},
$$
for fixed $a,b\in\F_q$ with $a\ne b$.

Theorem~\ref{t:n_m} below will cover a more general situation. To state the
result we fix some notations: let $a,b\in\F_q$, $a\ne b$, let $\beta\in\F_{q^m}$, and let 
$d=-1$ or $d=3$. In addition, if $d=3$ we assume $p\ne3$.
 
Let  
$$
L_{m,d,\beta}:=\F_{q^m}(x,y_a,y_b),\ y_a^q-y_a=x+a(\beta+x^d),\ y_b^q-y_b=x+b(\beta+x^d), 
$$
and for $u,v\ne0$ in the subfield $\F_q$ denote
$$
S_d^{(m)}(u,v)=\sum_{z}e(uz+vz^d),
$$
where $z$ runs over $\F_{q^m}^*$ or $\F_{q^m}$ according as $d$ equals $-1$ or $3$,
respectively. Moreover, we denote $S_d(u,v):=S_d^{(1)}(u,v)$.

\begin{thm}\label{t:n_m}
The number $N_{m,d}(\beta)$ of rational places of $L_{m,d,\beta}$ is
given by
$$
N_{m,d}(\beta)=q^m+1+\sum_{v\in\F_q^*}e(\beta v)\sum_{u}S_d^{(m)}(u,v),
$$
where $u$ runs over $\F_q^*$ or $\F_q$ according as $d$ equals $-1$ or $3$, respectively.
\end{thm}

\begin{proof}Let $\{u_1,\dots,u_r\}$ be a basis of $\F_q$ over $\F_p$. It
follows by Proposition 1.2 and by the proof of Proposition 1.1
in~\cite{Garcia91} that
$$
L_{m,d,\beta}=\F_{q^m}(x,y_{u_1},\dots,y_{u_r},z_{u_1},\dots,z_{u_r})
$$
with 
\begin{align*}
y_{u_i}^p-y_{u_i}&=u_i(x+a(\beta+x^d))=:f_i(x),\\ 
z_{u_i}^p-z_{u_i}&=u_i(x+b(\beta+x^d))=:g_i(x)
\end{align*}
for $i=1,\dots,r$. Let $L$ be the $\F_p$-subspace of $\F_{q^m}(x)$ spanned 
by the elements $f_1(x),\dots$, $f_r(x),g_1(x),\dots,g_r(x)$. Now, each
element $f(x)$ of $L$ is of the form
\begin{equation}\label{e:rep}
f(x)=(au+bv)\beta+(u+v)x+(au+bv)x^d,
\end{equation}
for some $u,v\in\F_q$.
 
If $f(x)=0$ then $u=v=0$ as $a\ne b$. It follows that the elements 
$f_1(x),\dots,g_r(x)$ are linearly independent over $\F_p$ and,
moreover, representation~(\ref{e:rep}) is unique.

Since the mapping $(u,v)\mapsto(u+v,au+bv)$ is linear and invertible, it is a
permutation of $\F_{q^2}\setminus\{\zero\}$, and therefore each non-zero $f(x)\in\F_q(x)$ is
of the form $f(x)=v\beta+ux+vx^d$ for unique $(u,v)\in\F_{q^2}\setminus\{\zero\}$.  

If $d=3$ Proposition~\ref{t:expsumplaces} implies
$$
N_{m,d}(\beta)=q^m+1+\sum_{(u,v)\in\F_q^2\setminus\{\zero\}}
\sum_{z\in\F_{q^m}}e(v\beta+uz+vz^d),
$$ 
and the claim follows now by noting that the inner sum equals zero if $v=0$.

Assume next that $d=-1$. Proposition~\ref{t:expsumplaces} now implies that
\begin{align}\label{a:dpoints}
N_{m,d}(\beta)=q^m+1&+\frac1{p-1}
\bigl(-(q^2-1)+\sum_{(u,v)\in\F_q^2\setminus\{\zero\}}|S_{u,v,\beta}|\bigr)\\
&+\sum_{(u,v)\in\F_q^2\setminus\{\zero\}}e(v\beta)
\sum_{z\in\F_{q^m}\setminus P_{u,v}}e(uz+vz^{-1}),\notag
\end{align}
where $S_{u,v,\beta}$ is the set of rational places of 
$$
\F_{q^m}(x,y),\ y^p-y=ux+v(\beta+x^{-1})
$$ 
lying above $P_{\infty}$, and $P_{u,v}$ is the set of poles of $f(x)=ux+v(\beta+x^{-1})$ in $\F_{q^m}$. By~\cite[Prop.III.7.8(c), Cor.III.3.8]{Stichtenoth} we know that 
\begin{equation}\label{e:infty}
|S_{u,v,\beta}|=
\begin{cases}
0&\text{ if $u=0$ and $\Tr(v\beta)\ne0$},\\
p&\text{ if $u=0$ and $\Tr(v\beta)=0$},\\
1&\text{ if }u\ne0.
\end{cases}
\end{equation}

Assume $\tr(\beta)=0$. Now, by~\eqref{a:dpoints} and \eqref{e:infty}, we get 
\begin{align*}
N_{m,d}(\beta)&=q^m+1+\frac1{p-1}\bigl(-(q^2-1)+(q-1)p+q^2-1-(q-1)\bigr)\\
&\hskip49pt+\sum_{v\in\F_q^*}\sum_{u\in\F_q^*}S_d^{(m)}(u,v)+
\sum_{v\in\F_q^*}\sum_{z\in\F_{q^m}^*}e(vz^{-1})\\
&=q^m+1+\sum_{v\in\F_q^*}\sum_{u\in\F_q^*}S_d^{(m)}(u,v).
\end{align*} 

Assume finally that $\tr(\beta)\ne0$. By~\eqref{a:dpoints} and
\eqref{e:infty} we now have  
\begin{align*}
N_{m,d}(\beta)&=q^m+1+\frac1{p-1}\bigl(-(q^2-1)+(q/p-1)p+q^2-1-(q-1)\bigr)\\
&\hskip50pt+\sum_{v\in\F_q^*}e(v\beta)\sum_{u\in\F_q^*}S_d^{(m)}(u,v)+
\sum_{v\in\F_q^*}e(v\beta)\sum_{z\in\F_{q^m}^*}e(vz^{-1})\\
&=q^m+1-1+\sum_{v\in\F_q^*}e(v\beta)\sum_{u\in\F_q^*}S_d^{(m)}(u,v)-\sum_{v\in\F_q^*}
\chi(v\tr(\beta)),
\end{align*}  
where $\chi$ is the canonical additive character of $\F_q$, and the claim
follows. 
\end{proof}

By the following result we see that in order to count $N_{m,d}(\beta)$ it is
enough to count $N_{m,d}(0)$ unless $\tr(\beta)\ne0$, $d=3$, and $r$ is even.
That case will be considered in the next section. From now on we use the
abbreviated notation $N_{m,d}:=N_{m,d}(0)$.
 
\begin{cor}\label{c:nperm}
Let $\beta_1,\beta_2\in\F_{q^m}$, and assume $\tr(\beta_1)=0$,
$\tr(\beta_2)\ne0$. Then $N_{m,d}(\beta_1)=N_{m,d}$. Moreover, if $d=-1$, or
$d=3$ and $r$ is odd, then
$$
N_{m,d}(\beta_2)=q^m+1-\frac{N_{m,d}-q^m-1}{q-1}.
$$
\end{cor}

\begin{proof}
By Theorem~\ref{t:n_m} it is clear that $N_{m,d}(\beta_1)=N_{m,d}$. Let
$\beta\in\F_{q^m}$. Since
$$
N_{m,3}(\beta)=q^m+1+\sum_{v\in\F_q^*}e(v\beta)\sum_{u\in\F_q^*}S_3^{(m)}(u,v)+
\underbrace{\sum_{v\in\F_q^*}S_3^{(m)}(0,v)}_{=0},
$$
we have in all cases
$$
N_{m,d}(\beta)=q^m+1+\sum_{v\in\F_q^*}e(v\beta)\sum_{u\in\F_q^*}\sum_{z}e(uz+vz^d).
$$
Now, by the substitution $z\mapsto u^{-1}z$ we get
$$
N_{m,d}(\beta)=q^m+1+\sum_{v\in\F_q^*}e(v\beta)\sum_{u\in\F_q^*}\sum_{z}e(z+vu^{-d}z^d),
$$
and since the map $u\mapsto vu^{-d}$ is a permutation of $\F_q$ we obtain
$$
N_{m,d}(\beta)=q^m+1+\Bigl(\underbrace{\sum_{v\in\F_q^*}e(v\beta)}_{=:S}\Bigr)
\Bigl(\sum_{u\in\F_q^*}S_d^{(m)}(1,u)\Bigr),
$$
where $S$ equals $q-1$ or $-1$ according as $\tr(\beta)$ is or is not zero,
respectively, and the claim follows now easily. 
\end{proof}
  
\begin{cor}If $q=2$ then
$$
N_{m,-1}-(2^m+1)=(-1)^{m-1}D_m(1,2)
=-\Bigl(\frac{-1+\sqrt{-7}}2\Bigr)^m-\Bigl(\frac{-1-\sqrt{-7}}2\Bigr)^m. 
$$
If $q=3$ then
\begin{multline*}
N_{m,-1}-(3^m+1)=2(-1)^{m-1}(D_m(-1,3)+D_m(2,3))\\
=-2\biggr(\Bigl(\frac{1+\sqrt{-11}}2\Bigr)^m+\Bigl(\frac{1-\sqrt{-11}}
2\Bigr)^m
+(-1+\sqrt{-2})^m+(-1-\sqrt{-2})^m\biggr). 
\end{multline*}
\end{cor}

\begin{proof}
Since $S_{-1}^{(m)}(u,v)=(-1)^{m-1}D_m(S_{-1}(u,v),q)$,
by~\cite[Thm.5.46]{Lidl97}, the claim follows now by Theorem~\ref{t:n_m} and
equation~\eqref{e:Dvalues}. 
\end{proof}

If $q$ is a power of two or three we can count $N_{m,-1}$ up to the evaluation
of Kronecker class numbers:  

\begin{prop}[{\cite[Lemma 7]{Moisio06c}}]\label{l:KD3}
Let $p=2$ or $p=3$, and let $q=p^r$ with $r\ge2$. Then 
$$
N_{m,-1}=q^m+1+(-1)^{m-1}(q-1)\sum_{t\in S_p}H(t^2-4q)D_m(t,q),
$$
where $H(d)$ is the Kronecker class number of $d$, and 
$$
S_p=\{t\in\Z\,:\,|t|< 2\sqrt q,\,t\equiv-1\,(e)\}
\text{ with } e=
\begin{cases}
4, &\text{if }p=2,\\ 
3, &\text{if }p=3.
\end{cases}
$$
\end{prop}

In the case where $q$ is unbounded power of two or three we are still able to
give the $N_{m,-1}$ provided that $m$ is relatively small:

\begin{prop}[{\cite[Theorem 13, Remark 6]{Moisio06c}}]\label{p:explicit3}
Let $q=p^r$ with $p=2$ or $p=3$ and $r\ge2$.
The number $N_{m,-1}$ of rational places of $L_{m,-1,0}$ is given by
\begin{center}
\begin{tabular}{c|c|c}
$m$ & $N_m'$ when $q=2^r$ & $N_m'$ when $q=3^r$\\ \hline
$1$ & $q$ & $q$\\
$2$ & $q^2$ & $q^2$\\
$3$ & $\pm q^2$ & $0$\\
$4$ & $0$ & $q^2$\\
$5$ & $(t_7\mp1)q^3$ & $\pm q^3$\\
$6$ & $\pm q^3$ & $(-1\pm1)q^3$\\
$7$ & $(t_9-t_7+1)q^4$ & $(u_9\mp1)q^4$\\
$8$ & $(1 \mp 1)q^4$ & $q^4-q+1$\\
$9$ & $(t_{11} - t_9-1)q^5$ & $(u_{11} - u_9)q^5$\\
$10$& $\tau_2q^2-q^5$ & $\tau_3q^2-q^5$
\end{tabular}
\end{center}
where $N_m'=(N_{m,-1}-q^m-1)/(q-1)+q-1$, $\pm=(-1)^r$, and
\begin{center}
\renewcommand{\arraystretch}{1.2}
\begin{tabular}{c|c}
$w^r+\overline{w}^r$ & $w$\\\hline
$t_7$ & $(1+\sqrt{-15})/4$\\
$t_9$ & $(-5+\sqrt{-39})/8$\\
$u_9$ & $(5+2\sqrt{-14})/9$\\
$u_{11}$ & $(-1+4\sqrt{-5})/9$\\
$\tau_2$ & $-3+\sqrt{-119}$\\
$\tau_3$ & $14+\sqrt{-1991}$\\ \hline
$t_{11}$ & $w_+^r+\overline{w_+}^r+w_-^r+\overline{w_-}^r$\\
$w_\pm$ & $(-3\pm\sqrt{505}+\sqrt{-510\mp 6\sqrt{505}})/32$
\end{tabular}
\end{center}
\end{prop}

\begin{rem}As noted in~\cite[Remarks 3 and 5]{Moisio06c} we may formulate
the differences of values of Ramanujan's tau-functions in the following way:
\begin{align*}
\tau(q)-2^{11}\tau(q/4)&=q^2(\mu_2^r+\bar\mu_2^r) \text{ with }
\mu_2=-3+\sqrt{-119} 
\text{ if } q=2^r,\\  
\tau(q)-3^{11}\tau(q/9)&=q^2(\mu_3^r+\bar\mu_3^r) \text{ with }
\mu_3=14+\sqrt{-1991}
\text{ if } q=3^r.
\end{align*}
\end{rem}

\section{The number of rational places of $L_{m,3,\beta}$}

In order to count $N_{m,3}(\beta)$ we need the following three results.

\begin{lemma}\label{lemma:Dickson}
Assume $p\ne3$. Then
$$
S_3^{(m)}(u,v)=(-1)^{m-1}D_m(S_3(u,v),q)\qquad\forall u,v\in\F_q,u\ne0.
$$
\end{lemma}

\begin{proof}By Weil's theorem (see e.g.~\cite[Thm.5.36]{Lidl97}) we know
that there exist complex numbers $\omega$ and $\nu$ satisfying $|\omega|=|\nu|=\sqrt q$
and 
$$
S_3^{(m)}(u,v)=-(\omega^m+\nu^{m}).
$$
Obviously $S_3^{(m)}(u,v)$ is real and it follows that $\nu=\bar\omega$. Now,
by~(\ref{e:Dfunct})
$$
D_m(\omega+q/\omega,q)=\omega^m+q^m/\omega^m=\omega^m+\bar\omega^{m},
$$
and therefore
$$
S_3^{(m)}(u,v)=-D_m(\omega+\bar\omega,q)=-D_m(-S_3(u,v),q).
$$
\end{proof}

\begin{prop}\label{prop:values} 
Let $q=2^r$. The value distribution of $S_3(u,v)$ is given by
\begin{center} 
\begin{tabular}{c | c | c}
$S_3(u,v)$ & $\#$  in case $2\mid r$ & $\#$  in case $2\nmid r$ \\ \hline
$q$ & $1$ & $1$\\
$-2\sqrt{q}$ & $\frac{q-1}{24}(q-2\sqrt{q})$ & $0$\\
$-\sqrt{2q}$ & $0$ & $\frac{q-1}4(q-\sqrt{2q})$\\
$-\sqrt{q}$ & $\frac{q-1}{3}(q-\sqrt{q})$ & $0$\\
$0$ & $\frac{q-1}4 q+(q-1)$ & $\frac{q-1}2q +(q-1)$ \\
$\sqrt{q}$ & $\frac{q-1}{3}(q+\sqrt{q})$ & $0$\\
$\sqrt{2q}$ & $0$ & $\frac{q-1}4(q+\sqrt{2q})$\\
$2\sqrt{q}$ & $\frac{q-1}{24}(q+2\sqrt{q})$ & $0$
\end{tabular}
\end{center}
where $S_3(u,v)$ is attained $\#$ times as $(u,v)$ varies over $\F_q^2$. 
\end{prop}

\begin{proof}Let $\gamma$ be a primitive element of $\F_q$, and let 
$$
c(u,v)=(\Tr(u+v),\Tr(u\gamma+v\gamma^3),\dots,\Tr(u\gamma^{q-2}+v\gamma^{3(q-2)})),
$$ 
be a codeword in the dual $B^{\perp}$ of the binary two-error-correcting
BCH code of length 
$q-1$, and let $w(c(u,v))$ denote the Hamming weight of $c(u,v)$.

The claim follows now by the weight distribution of $B^{\perp}$ (see
e.g~\cite{Schoof95}) and by the following two facts which are easy to verify:
\begin{enumerate}
\item Map $\psi:(\F_q^2,+)\rightarrow B^{\perp}, (u,v)\mapsto c(u,v)$ is a
group isomorphism.
\item $S(u,v)=q-2w(c(u,v))$.
\end{enumerate} 
\end{proof}

We rephrase a result by Carlitz~\cite{Carlitz79} in the following form.

\begin{prop}\label{prop:Carlitz}
Let $q=2^r$, $\gamma$ be a primitive element $\F_q$, and let
$v\in\F_q^*$. If $r$ is even, then 
\[S_3(u,v)\in\begin{cases}
\{0,\pm 2\sqrt q\} & \text{if }v\in\langle \gamma^3\rangle,\\
\{\pm \sqrt q\} & \text{if }v\not\in\langle \gamma^3\rangle,
\end{cases}\]
and each value is attained at least once as $u$ varies over $\F_q$.
Moreover,
\[S_3(0,v)=\begin{cases}
(-1)^{r/2+1}2\sqrt q & \text{if }v\in\langle \gamma^3\rangle,\\
(-1)^{r/2}\sqrt q & \text{if }v\not\in\langle \gamma^3\rangle.
\end{cases}\]
\end{prop}

Now we have all the tools in order to establish the main result of this section: 
 
\begin{thm}\label{t:expd3}
Let $q=2^r$, let $\beta\in\F_{q^m}$, and let $c=\tr(\beta)$. If $r$ is odd,
the number $N_{m,3}$ of rational places of $L_{m,3,0}$ is given by
\begin{center}
\begin{tabular}{c | c}
$m \mod 8$ & $N_{m,3}-(q^m+1)$  \\ \hline
$0$ & $-2(q-1)q^{\frac{m+2}2}$\\
$\pm1$ & $(q-1)q^{\frac{m+1}2}$\\
$\pm2$ & $(q-1)q^{\frac{m+2}2}$\\
$\pm3$ & $-(q-1)q^{\frac{m+1}2}$\\
$4$ & $0$\\
\end{tabular}
\end{center}
\vskip3pt
If $r$ is even, $N_{m,3}$ is given by
\begin{center}
\begin{tabular}{c | c}
$m \mod 12$ & $N_{m,3}-(q^m+1)$  \\ \hline
$0$ & $-2(q-1)q^{\frac{m+2}2}$\\
$\pm1,\pm5$ & $(q-1)q^{\frac{m+1}2}$\\
$\pm2$ & $(q-1)q^{\frac{m+2}2}$\\
$\pm3$ & $-(q-1)q^{\frac{m+1}2}$\\
$\pm4$ & $0$\\
$6$ & $-(q-1)q^{\frac{m+2}2}$
\end{tabular}
\end{center}
\vskip3pt
If $r$ is even and $c\ne0$, the number $N_{m,3}(\beta)$ of rational places
of $L_{m,3,\beta}$ is given by
\begin{center}
\begin{tabular}{c|c|c}
$m \mod 12$ & $N_{m,3}(\beta)-(q^m+1)$, $c\in\langle\gamma^3\rangle$ 
& $N_{m,3}(\beta)-(q^m+1)$, $c\not\in\langle\gamma^3\rangle$\\ \hline
$0$ & $2q^\frac{m+2}{2}$ & $2q^\frac{m+2}{2}$\\
$\pm1,\pm5$ & $-q^\frac{m+1}{2}$ & $-q^\frac{m+1}{2}$\\
$\pm2$ & $-q^\frac{m+2}{2}$ & $-q^\frac{m+2}{2}$\\
$\pm3$ & $(1-(-1)^s 2\sqrt{q})q^\frac{m+1}{2}$ &
  $(1+(-1)^s \sqrt{q})q^\frac{m+1}{2}$  \\
$\pm4$ & $(-1)^s2q^\frac{m+3}{2}$ & $(-1)^{s+1}q^\frac{m+3}{2}$\\
$6$ &  $(1-(-1)^s 2\sqrt{q})q^\frac{m+2}{2}$ & 
$(1+(-1)^s \sqrt{q})q^\frac{m+2}{2}$
\end{tabular}
\end{center}
where $s=r/2$.
\end{thm}

\begin{proof}Assume first that $r$ is odd. Now, By Theorem~\ref{t:n_m},
Lemma~\ref{lemma:Dickson}, and Proposition~\ref{prop:values} we obtain
\begin{align*}
N_{m,3}(\beta)-(q^m+1)
=&\sum_{v\in\F_q^*}\sum_{u\in\F_q}S_3^{(m)}(u,v)\\
=&(-1)^{m-1}\Big(\textstyle
\frac{q-1}{4}(q-\sqrt{2q})D_m(-\sqrt{2q},q)\\
&+\textstyle\frac{q-1}{4}(q+\sqrt{2q})D_m(\sqrt{2q},q)
+\textstyle\frac{q-1}{2}qD_m(0,q)\Big).
\end{align*}
We note that above the $q-1$ zeros of $S_3(u,v)$ corresponding the pairs 
$(u,0)$, $u\ne 0$, are excluded. By \eqref{e:Dvalues} we see that 
$D_m(\pm \sqrt{2q},q)=2(\pm q)^m\cos(\frac{m\pi}4)$, $D_m(0,q)=0$ if $m$ is odd, and 
$D_m(0,q)=2(-q)^{m/2}$ if $m$ is even. Thus
$$
N_{m,3}(\beta)-(q^m+1)=
\begin{cases}
\sqrt2(q-1)q^{\frac{m+1}2}\cos(\frac{m\pi}4)&\text{ if }2\nmid m,\\
-(q-1)q^{\frac{m+1}2}\big(\cos(\frac{m\pi}4)+(-1)^{\frac m2}\big)&\text{ if }2\mid m,  
\end{cases}
$$
and the claim follows.
 
Assume next that $r$ is even. Let $\chi$ be the canonical additive character of $\F_q$. By
Theorem~\ref{t:n_m}
\begin{align*}
S:=&N_{m,3}(\beta)-(q^m+1)\\
=&\sum_{v\in\F_q^*}\chi(cv)\sum_{u\in\F_q}S_3^{(m)}(u,v)\\
=&\sum_{i=0}^2\sum_{v\in\gamma^i\langle\gamma^3\rangle}
\chi(cv)\sum_{u\in\F_q}S_3^{(m)}(u,v)\\
=&\sum_{i=0}^2\sum_{j=0}^{(q-4)/3}
\chi(c\gamma^{i+3j})\underbrace{\sum_{u\in\F_q}\sum_{x\in\F_{q^m}}
e(ux+\gamma^i(\gamma^jx)^3)}_{=:S^*}.
\end{align*}

By the substitution $x\mapsto \gamma^{-j}x$ we have
\[S^*=\sum_{u\in\F_q}\sum_{x\in\F_{q^m}}
e(u\gamma^{-j}x+\gamma^ix^3)\overset{u\mapsto \gamma^ju}{=}
\sum_{u\in\F_q}S_3^{(m)}(u,\gamma^i),\]
and therefore
\[S=\sum_{j=0}^{(q-4)/3}\chi(c\gamma^{3j})\sum_{u\in\F_q}S_3^{(m)}(u,1)
+\sum_{i=1}^2\sum_{j=0}^{(q-4)/3}
\chi(c\gamma^{i+3j})\sum_{u\in\F_q}S_3^{(m)}(u,\gamma^i).\]

We observe that
\begin{align} \label{eq:subst}
\sum_{u\in\F_q}S_3^{(m)}(u,\gamma^2)\overset{x\mapsto x^2}{=}&
\sum_{u\in\F_q}\sum_{x\in\F_{q^m}}e(ux^2+\gamma^2(x^3)^2)\\
\overset{u\mapsto u^2}{=}&
\sum_{u\in\F_q}\sum_{x\in\F_{q^m}}e((ux+\gamma x^3)^2)=
\sum_{u\in\F_q}S_3^{(m)}(u,\gamma), \notag
\end{align}
and we now get
\[S=\underbrace{\sum_{j=0}^{(q-4)/3}\chi(c\gamma^{3j})}_{S_1} 
\underbrace{\sum_{u\in\F_q}S_3^{(m)}(u,1)}_{S_2}
+\underbrace{\sum_{i=1}^2\sum_{j=0}^{(q-4)/3}\chi(c\gamma^{i+3j})}_{S_3}
\underbrace{\sum_{u\in\F_q}S_3^{(m)}(u,\gamma)}_{S_4}.\]

Consider each sum $S_i$ above separately. Clearly, $S_1=(S_3(0,c)-1)/3$, and
by Proposition~\ref{prop:Carlitz}
\begin{align*}
S_3=&\frac{1}{3}(S_3(0,c\gamma)-1+S_3(0,c\gamma^2)-1)\\
=&
\begin{cases}
\frac23(q-1) &\text{if }c=0,\\
\frac{2}{3}((-1)^sq^{1/2}-1) &\text{if }c\in\langle\gamma^3\rangle,\\
\frac{1}{3}((-1)^{s+1}2q^{1/2}+(-1)^{s}q^{1/2}-2) &\text{if
}c\not\in\langle\gamma^3\rangle.\end{cases}
\end{align*}

We apply the argument used with the sum $S^*$ above to the opposite
direction to get
\[S_2=\frac{3}{q-1}\sum_{j=0}^{(q-4)/3}
\sum_{u\in\F_q}S_3^{(m)}(u,\gamma^{3j})
=\frac{3}{q-1}\sum_{v\in\langle\gamma^3\rangle}
\sum_{u\in\F_q}S_3^{(m)}(u,v).\]
By Proposition \ref{prop:Carlitz} we know that the sum $S_3(u,v)$ gets
exactly the values $0,\pm 2\sqrt{q}$ when $v\in\langle\gamma^3\rangle$. By
Lemma \ref{lemma:Dickson} and Proposition \ref{prop:values} we obtain
\begin{align*}
\frac{q-1}{3}S_2=&(-1)^{m-1}\Big(\textstyle
\frac{q-1}{24}(q-2\sqrt{q})D_m(-2\sqrt{q},q)\\
&+\textstyle\frac{q-1}{24}(q+2\sqrt{q})D_m(2\sqrt{q},q)
+\textstyle\frac{q-1}{4}qD_m(0,q)\Big).
\end{align*}
By \eqref{e:Dvalues} we see that $D_m(\pm 2\sqrt{q},q)= 2(\pm\sqrt{q})^m$, 
$D_m(0,q)=0$ if $m$ is odd, and $D_m(0,q)=2(-q)^{m/2}$ if $m$ is even. 
All in all,
\[S_2=\begin{cases}
-\frac{1}{2}q^{\frac{m}{2}+1}-\frac{3}{2}q(-q)^\frac{m}{2} & \text{if } 2\mid
m,\\
q^\frac{m+1}{2} & \text{if } 2\nmid m.\end{cases}\]

Consider finally $S_4$. Now
\[S_4=\frac{1}{2}\sum_{i=1}^2 \sum_{u\in\F_q} S_3^{(m)}(u,\gamma^i)=
\frac{3}{2(q-1)}\sum_{v\in\F_q^*\setminus\langle\gamma^3\rangle}
\sum_{u\in\F_q}S_3^{(m)}(u,v)\]
and by Proposition \ref{prop:Carlitz} the value set of $S_3(u,v)$ is 
$\{\pm\sqrt{q}\}$. Again, by Lemma \ref{lemma:Dickson} and Proposition
\ref{prop:values}, we get
\begin{align*}
\textstyle\frac{2(q-1)}{3}S_4=&(-1)^{m-1}\Big(\textstyle
\frac{q-1}{3}(q-\sqrt{q})D_m(-\sqrt{q},q)
+\textstyle\frac{q-1}{3}(q+\sqrt{q})D_m(\sqrt{q},q)\Big).
\end{align*}
By \eqref{e:Dvalues} $D_m(\pm\sqrt{q},q)$ equals $ 2q^{m/2}$, $\pm
q^{m/2}$, $- q^{m/2}$, and $\mp 2q^{m/2}$ when $m\equiv 0$, $\pm1$,
$\pm 2$, and $3\pmod{6}$, respectively. Hence, we get
\[S_4=\begin{cases}
-2q^{\frac{m}{2}+1} & \text{if } m\equiv 0\pmod{6},\\
q^\frac{m+1}{2} & \text{if } m\equiv \pm 1\pmod{6},\\    
q^{\frac{m}{2}+1} & \text{if } m\equiv \pm 2\pmod{6},\\
-2q^\frac{m+1}{2} & \text{if } m\equiv 3\pmod{6}.\end{cases}
\]
By collecting all the calculations we obtain the claimed result.
\end{proof}

\section{Enumeration of irreducible polynomials with prescribed coefficients}

In this section we calculate the number of irreducible polynomials over
$\F_q$ in the cases (i) and (ii) of the Introduction. The method we use here
is a modification of the method introduced in~\cite{Niederreiter90} (see also
\cite{Cohen00,Chou01}). Roughly speaking, the method involves two steps:
first, count the number of all the elements of $\F_{q^m}$ with prescribed
traces, and then, by using M\"obius inversion, count the number of elements of
degree $m$ with prescribed coefficients. From now on we assume that $p=2$ or
$p=3$.

\subsection{Elements of degree $m$ with prescribed traces}

Let $d=3$ or $d=-1$ and employ the convention $0^{-1}=0$.

\begin{defn}
For $c\in\F_q$ define
$$
H_{c,d}(m)=|\{z\in\F_{q^m}\ |\ \tr(z)=0,\ \tr(z^d)=c\}|.
$$
\end{defn}

\begin{lemma}\label{l:n_elem}
Let $\alpha\in\F_{q^m}$ satisfying $\tr(\alpha)=1$. Then
$$
H_{c,d}(m)=\frac1{q^2}(N_{m,d}(-\alpha c)-1+\epsilon),
$$
where 
$$
\epsilon=
\begin{cases}
0&\text{ if }d=3,\\
1-q&\text{ if $d=-1$ and $c\ne0$},\\
(q-1)^2&\text{ if $d=-1$ and $c=0$}.
\end{cases}
$$
\end{lemma}

\begin{proof}Let $\chi$ be the canonical additive character of $\F_q$. Now 
\begin{align*}\label{eq:vali}
H_{c,d}(m)=&\sum_{z\in\F_{q^m}}
\Big(\frac{1}{q}\sum_{u\in\F_q}\chi(\tr(z)u)\Big) 
\Big(\frac{1}{q}\sum_{v\in\F_q}\chi(\tr(z^d-\alpha c)v)\Big)\\
=& \frac{1}{q^2}\sum_{u,v\in\F_q} \sum_{z\in\F_{q^m}}
e(uz+vz^d-\alpha cv)\\
=&\frac1{q^2}\Big(q^m+\sum_{v\in\F_q^*}e(-\alpha
cv)\sum_{u\in\F_q}\sum_{z\in\F_{q^m}}
e(uz+vz^d)\Big)\\
=&\frac1{q^2}\Big(q^m+\sum_{v\in\F_q^*}e(-\alpha c
v)\sum_{u}S^{(m)}_d(u,v)+\epsilon\Big),
\end{align*}
where the last equality follows by the definition of $S_{d}(u,v)$. The claim
follows now by Theorem~\ref{t:n_m}.
\end{proof}

Next we shall count the number of elements $z$ in $\F_{q^m}$ of degree $m$
satisfying $\tr(z)=0$ and $\tr(z^d)=c$.
 
\begin{defn} For $c\in\F_q$ define
$$
G_{c,d}(m)=|\{z\in\F_{q^m}\ |\ \tr(z)=0,\
\tr(z^d)=c,\ z\not\in\F_{q^n}^*\text{ if }n<m\}|.
$$
\end{defn}

We shall need the following well-known formula for the number of all irreducible
polynomials (see e.g. \cite[Thm.3.25]{Lidl97}):

\begin{prop}
The number of monic irreducible polynomials in $\F_q[x]$ of degree $m$
is given by $I(m)/m$, where
\[I(m)=\sum_{t|m}\mu(t)q^\frac{m}{t}.\]
\end{prop}

Let $n$ be a positive factor of $m$, and let $\tr_n:\F_{q^n}\to\F_q$ denote
the relative trace function. Clearly, for every $z\in\F_{q^n}$ we have
$\tr(z)=\frac{m}{n}\tr_n(z)$ and therefore
$$
\tr(z)=\tr(z^d)=0 \text{ iff }\big[
p\mid{\textstyle\frac{m}{n}} \text{ or } \big(p\nmid\textstyle\frac{m}{n}
\text{ and }\tr_n(z)=\tr_n(z^d)=0\big) \big].
$$
\vskip5pt
Let $m=p^k s$ such that $p\nmid s$. Now
\begin{align*}
H_{0,d}(m)=&H_{0,d}(p^ks)=\sum_{n|m,\ p|\frac{m}{n}}I(n)+
\sum_{n|m,\ p\,\nmid\frac{m}{n}}G_{0,d}(n)\\
=&\sum_{t|s}\sum_{i=0}^{k-1}I(p^i t)+\sum_{t|s}G_{0,d}(p^kt)
=\sum_{t|s}\big(S(p^kt)+G_{0,d}(p^kt)\big),
\end{align*}
where $S(p^kt)=\sum_{i=0}^{k-1}I(p^i t)$. By M\"obius inversion, see e.g.
\cite[Thm.3.24]{Lidl97}, we get
\[S(p^ks)+G_{0,d}(p^ks)=\sum_{t|s}\mu\big(\frac{s}{t}\big)H_{0,d}(p^kt).\]
By Lemma~\ref{l:n_elem} we now get the following theorem.

\begin{thm}\label{thm:G0}
Let $m=p^ks$ with coprime $p$ and $s$. Then
\[G_{0,d}(m)=\sum_{t|s}\mu\big(\frac{s}{t}\big)H_{0,d}(p^kt)-S(m),\]
where
\[H_{0,d}(n)=\frac{1}{q^2}(N_{n,d}(0)-1+\epsilon)\ \text{ and }\ S(m)=
\sum_{i=0}^{k-1}I(p^is)\]
with
$$
\epsilon=
\begin{cases}
0&\text{ if }d=3,\\
(q-1)^2&\text{ if }d=-1.
\end{cases}
$$
\end{thm}

Assume next that $c\ne0$, and let $n$ be a positive factor of $m$. Now, for
each $z\in\F_{q^n}$, we have that 
$$
\tr(z)=0\text{ and }\tr(z^d)=c\text{ iff }\tr_n(z)=0\text{ and}\tr_n(z^d)
=\textstyle\frac nm c.
$$
Since $n/m=1$ or $n/m=\pm1$ according as $p$ equals $2$ or $3$, respectively, we see
that
$$
G_{\frac{nc}m,d}(n)=G_{c,d}(n),
$$ 
and therefore
\[H_{c,d}(m)=\sum_{t|s}G_{c,d}(p^kt).\]

Now, by M\"obius inversion and Lemma~\ref{l:n_elem} we get 
\begin{thm}\label{thm:Gc}
Let $m=p^ks$ with $p$ and $s$ coprime, let $c\in\F_q^*$, and let
$\alpha\in\F_{q^m}$ satisfying $tr(\alpha)=1$. Then 
\[G_{c,d}(m)=\sum_{t|s}\mu\big(\frac{s}{t}\big)H_{c,d}(p^kt),\]
where
\[H_{c,d}(n)=\frac{1}{q^2}\Big(N_{n,d}(-\alpha c)-1+\epsilon\Big)\]
with 
$$
\epsilon=
\begin{cases}
0&\text{ if }d=3,\\
1-q&\text{ if }d=-1.
\end{cases}
$$
\end{thm}

\subsection{Irreducible polynomials of degree $m$ with prescribed traces}

\begin{lemma}\label{l:p_1}
Let $q=p^r$ with $p=2$ or $p=3$, and let $c\in\F_q$. The number of irreducible
polynomials $p(x)=x^m+a_{m-1}x^{m-1}+\cdots+a_1x+a_0$ in $\F_q[x]$ with
$a_{m-1}=c$ and $a_1=0$ equals $G_{c,-1}(m)/m$. 
\end{lemma}

\begin{proof}
If $p(x)$ is irreducible, then $a_{m-1}=\tr(z)$, where $z$ is any
of the $m$ distinct roots of $p(x)$ in $\F_{q^m}$. Moreover, since
$a_0^{-1}x^mp(x^{-1})$ is monic and irreducible we get $\tr(z^{-1})=a_1/a_0$.
Hence, the number of irreducible $p(x)$ with $a_{m-1}=c$ and $a_1=0$ equals
$G_{c,-1}^{\prime}(m)/m$, where $G_{c,-1}^{\prime}(m)$ is the number of
elements $z$ of degree $m$ over $\F_q$ in $\F_{q^m}$ satisfying $\tr(z)=c$
and $\tr(z^{-1})=0$. But clearly $G_{c,-1}^{\prime}(m)=G_{c,-1}(m)$. 
\end{proof}

\begin{rem}
By the preceding proof it is clear that the number of irreducible polynomials
$x^m+a_{m-1}x^{m-1}+\dots+a_1x+a_0$ in $\F_q[x]$ such that $a_{m-1}=0$ and $a_1/a_0=c$ also 
equals $G_{c,-1}(m)/m$. 
\end{rem}

\begin{lemma}Let $q=2^r$, and let $c\in\F_q$. The number of irreducible polynomials $p(x)=x^m+a_{m-1}x^{m-1}+\cdots+a_1x+a_0$ in $\F_q[x]$ with $a_{m-1}=0$ and $a_{m-3}=c$ equals 
$G_{c,3}(m)/m$. 
\end{lemma}

\begin{proof}\label{l:p_3}
Let $z=z_1,\dots,z_m$ be the roots of an irreducible polynomial $p(x)$ in
$\F_{q^m}$. If $m\ge3$, we have, by Newton's formula (see~\cite[Thm.1.75]{Lidl97}), that
$$
s_3+s_2a_{m-1}+s_1a_{m-2}+a_{m-3}=0,
$$
where $s_k=\sum_{i=1}^mz_i^k=\tr(z^k)$. Since $s_1=a_1=0$, we get $a_{m-3}=\tr(z^3)$, and the claim follows.    
\end{proof}

Lemma~\ref{l:p_1}, Theorems~\ref{thm:G0} and~\ref{thm:Gc},
Corollary~\ref{c:nperm}, and Proposition~\ref{p:explicit3} give, with help of
{\em Mathematica}, the following two 
corollaries:
 
\begin{cor}
Let $q=p^r$ with $p=2$ or $p=3$. The number of irreducible
polynomials $x^m+a_{m-1}x^{m-1}+\dots+a_1x+a_0$ in $\F_q[x]$ with
$a_{m-1}=a_1=0$ equals $G_{0,-1}(m)/m$, where 
\begin{center}
\renewcommand{\arraystretch}{1}
\begin{tabular}{c|c}
$m$ & $G_{0,-1}(m)$ with $q=2^r$\\ \hline
$1$ & $1$\\
$2$ & $0$\\
$3$ & $(1\pm 1)(q-1)$\\
$4$ & $0$\\
$5$ & $q^3+(t_7\mp 1)q(q-1)-1$\\
$6$ & $(q-1)(q^3\pm q)$\\
$7$ & $q^5+q^2(q-1)(t_9-t_7+1)-1$\\
$8$ & $q^6-q^4+(1\mp 1)q^2(q-1)$\\
$9$ & $q^7+(q-1)(q^3(t_{11}-t_9-1)-1\mp 1)-1$\\
$10$ & $q^8-q^5-q^4+q^3+(q-1)\tau_2$\\
\end{tabular}
\vskip5pt
\begin{tabular}{c|c}
$m$ & $G_{0,-1}(m)$ with $q=3^r$ \\ \hline
$1$ & $0$\\
$2$ & $q-1$\\
$3$ & $0$\\
$4$ & $q^2-1$\\
$5$ & $q^3\pm q(q-1)-1$\\
$6$ & $q(q-1)(q^2+q-1\pm 1)$\\
$7$ & $q^5+q^2(q-1)(u_9\mp 1)-1$\\
$8$ & $q^6+q^3-2q^2-q+1$\\
$9$ & $q^7+q^3(q-1)(u_{11}-u_9)-q^3$\\
$10$ & $q^8-q^4-(q-1)(\pm q+1)+(q-1)\tau_3$\\
\end{tabular}
\end{center}
\end{cor}

\begin{cor}
Let $q=p^r$ with $p=2$ or $p=3$. The number of monic irreducible
polynomials $x^m+a_{m-1}x^{m-1}+\dots+a_1x+a_0$ in $\F_q[x]$ with
$a_{m-1}=c\ne 0$ and $a_1=0$ equals $G_{c,-1}(m)/m$, where
\begin{center}
\renewcommand{\arraystretch}{1}
\begin{tabular}{c|c|c}
$m$ & $G_{c,-1}(m)$ with $q=2^r$ & $G_{c,-1}(m)$ with $q=3^r$ \\ \hline
$1$ & $0$ & $0$\\
$2$ & $0$ & $0$\\
$3$ & $q\mp 1$ & $q$\\
$4$ & $q^2$ & $q^2-1$\\
$5$ & $q^3-q(t_7\mp 1)$ & $q^3\mp q$\\
$6$ & $q^4 \mp q$ & $q^4\mp q$\\
$7$ & $q^5-q^2(t_9-t_7+1)$ & $q^5-q^2(u_9\mp 1)$\\
$8$ & $q^6+(-1\pm 1)q^2$ & $q^6-2q^2+1$\\
$9$ & $q^7-q^3(t_{11}-t_9-1)-q\pm 1$ & $q^7-q^3(u_{11}-u_9)$\\
$10$ & $q^8+q^3-\tau_2$ &
$q^8\pm q-\tau_3$\\
\end{tabular}
\end{center}
\end{cor}

Lemma~\ref{l:p_3}, Theorems~\ref{thm:G0} and~\ref{thm:Gc},
Corollary~\ref{c:nperm}, and Theorem~\ref{t:expd3} give, with help of {\em
Mathematica}, the following  

\begin{cor}
Let $q=2^r$, and let $c\in\F_q$. The number of monic irreducible
polynomials $x^m+a_{m-1}x^{m-1}+\dots+a_1x+a_0$ in $\F_q[x]$ with
$a_{m-1}=0$ and $a_{m-3}=c$ equals $G_{c,3}(m)/m$, where
\begin{center}
\tiny
\renewcommand{\arraystretch}{1}
\begin{tabular}{c|c|c|c}
$m$ & $G_{c,3}(m)$, if $2\mid r$ and $c=0$ & $G_{c,3}(m)$, if $2\nmid r$ and  $c=0$ & $G_{c,3}(m)$, if $2\nmid r$ and $c\ne0$  \\ \hline
$1$ & $1$ & $1$ & $0$\\
$2$ & $0$ & $0$ & $0$\\
$3$ & $0$ & $0$ & $q+1$\\
$4$ & $0$ & $0$ & $q^2$\\
$5$ & $q^3+q^2-q-1$ & $q^3-q^2+q-1$ & $q^3+q$\\
$6$ & $q^4-2q^3+q^2$ & $q^4-q^2$ & $q^4-q^2$\\
$7$ & $q^5+q^3-q^2-1$ & $q^5+q^3-q^2-1$ & $q^5-q^2$\\
$8$ & $q^6-q^4$ & $q^6-3q^4+2q^3$ & $q^6+2q^3$\\
$9$ & $q^7-q^4+q^3-1$ & $q^7+q^4-q^3-1$ & $q^7-q^3-q-1$\\
$10$ & $q^8-q^4$ & $q^8-q^4$ & $q^8-q^4$\\
$11$ & $q^9+q^5-q^4-1$ & $q^9-q^5+q^4-1$ & $q^9+q^4$\\
$12$ & $q^{10}-3q^6+2q^5$ & $q^{10}-q^6$ & $q^{10}-q^2$\\
$13$ & $q^{11}+q^6-q^5-1$ & $q^{11}-q^6+q^5-1$ & $q^{11}+q^5$\\
$14$ & $q^{12}-q^6$ & $q^{12}-q^6$ & $q^{12}-q^6$\\
$15$ & $q^{13}-q^7+q^6-q^3-q^2+q$ & $q^{13}+q^7-q^6-q^3+q^2-q$ &
$q^{13}-q^6-q^3-2q-1$\\
$16$ & $q^{14}-q^8$ & $q^{14}-3q^8+2q^7$ & $q^{14}+2q^7$\\
$17$ &$q^{15}+q^8-q^7-1$ & $q^{15}+q^8-q^7-1$ & $q^{15}-q^7$\\
$18$ & $q^{16}-2q^9+q^8-q^4+2q^3-q^2$ & $q^{16}-q^8-q^4+q^2$ &
$q^{16}-q^8-q^4+q^2$\\
$19$ & $q^{17}+q^9-q^8-1$ & $q^{17}-q^9+q^8-1$ & $q^{17}+q^8$\\
$20$ & $q^{18}-q^{10}$ & $q^{18}-q^{10}$ & $q^{18}-q^2$\\
$21$ & $q^{19}-q^{10}+q^9-q^5-q^3+q^2$  & $q^{19}-q^{10}+q^9-q^5-q^3+q^2$ & 
       $q^{19}+q^9-q^5+q^2-q-1$\\
$22$ & $q^{20}-q^{10}$  & $q^{20}-q^{10}$ & $q^{20}-q^{10}$ \\
$23$ & $q^{21}+q^{11}-q^{10}-1$  & $q^{21}+q^{11}-q^{10}-1$ &
       $q^{21}-q^{10}$\\
$24$ & $q^{22}-3q^{12}+2q^{11}-q^6+q^4$ &
       $q^{22}-3q^{12}+2q^{11}-q^6+3q^4-2q^3$ &
       $q^{22}+2q^{11}-q^6-2q^3$ \\
$25$ & $q^{23}+q^{12}-q^{11}-q^3-q^2+q$ & $q^{23}+q^{12}-q^{11}-q^3+q^2-q$ &
       $q^{23}-q^{11}-q^3-q$\\
$26$ & $q^{24}-q^{12}$ & $q^{24}-q^{12}$ & $q^{24}-q^{12}$ \\
$27$ & $q^{25}-q^{13}+q^{12}-q^7+q^4-q^3$ &
       $q^{25}-q^{13}+q^{12}-q^7-q^4+q^3$ &
       $q^{25}+q^{12}-q^7+q^3$ \\
$28$ & $q^{26}-q^{14}$ & $q^{26}-q^{14}$ & $q^{26}-q^2$\\
$29$ & $q^{27}+q^{14}-q^{13}-1$  & $q^{27}-q^{14}+q^{13}-1$ &
       $q^{27}+q^{13}$ \\
$30$ & $q^{28}-2q^{15}+q^{14}-q^8+2q^3-q^2$ & $q^{28}-q^{14}-q^8+q^2$ &
       $q^{28}-q^{14}-q^8+q^2$
\end{tabular}
\end{center}
moreover, if $r$ is even and $c\ne0$ then
\begin{center}
\tiny
\renewcommand{\arraystretch}{1}
\begin{tabular}{c|c|c}
$m$ & $G_{c,3}(m)$, if $c\in<\gamma^3>$ & $G_{c,3}(m)$, if $c\not\in<\gamma^3>$\\ \hline
$1$ & $0$ & $0$ \\
$2$ & $0$ & $0$ \\
$3$ & $q\mp2q^\frac{1}{2}+1$ & $q\pm q^\frac{1}{2}+1$\\
$4$ & $q^2 \pm 2 q^\frac{3}{2}$ & $q^2 \mp q^\frac{3}{2}$ \\
$5$ & $q^3-q$ & $q^3-q$\\
$6$ & $q^4\mp 2 q^\frac{5}{2}+q^2$ & $q^4\pm q^\frac{5}{2}+q^2$\\
$7$ & $q^5-q^2$ & $q^5-q^2$\\
$8$ & $q^6\pm 2q^\frac{7}{2}$ & $q^6\mp q^\frac{7}{2}$\\
$9$ & $q^7 \mp 2q^\frac{7}{2}+q^3-q\pm 2 q^\frac{1}{2}-1$ & 
      $q^7 \pm q^\frac{7}{2}+q^3-q\mp q^\frac{1}{2}-1$\\
$10$ & $q^8-q^4$ & $q^8-q^4$\\
$11$ & $q^9-q^4$ & $q^9-q^4$\\
$12$ & $q^{10}+2q^5-q^2\mp 2q^\frac{3}{2}$ & 
       $q^{10}+2q^5-q^2\pm q^\frac{3}{2}$\\
$13$ & $q^{11}-q^5$ & $q^{11}-q^5$\\
$14$ & $q^{12}-q^6$ & $q^{12}-q^6$\\
$15$ & $q^{13}\mp 2 q^\frac{13}{2}+q^6-q^3\pm 2 q^\frac{1}{2}-1$ & 
       $q^{13}\pm q^\frac{13}{2}+q^6-q^3\mp q^\frac{1}{2}-1$\\
$16$ & $q^{14}\pm 2q^\frac{15}{2}$ & $q^{14}\mp q^\frac{15}{2}$\\
$17$ & $q^{15}-q^7$ & $q^{15}-q^7$\\
$18$ & $q^{16} \mp 2q^\frac{17}{2}+q^8-q^4\pm 2q^\frac{5}{2}-q^2$ &
       $q^{16} \pm q^\frac{17}{2}+q^8-q^4\mp q^\frac{5}{2}-q^2$\\
$19$ & $q^{17}-q^8$ & $q^{17}-q^8$\\
$20$ & $q^{18} \pm 2 q^\frac{19}{2}-q^2 \mp 2q^\frac{3}{2}$ & 
       $q^{18} \mp q^\frac{19}{2}-q^2 \pm q^\frac{3}{2}$\\
$21$ & $q^{19}\mp 2 q^\frac{19}{2}+q^9-q^5+q^2-q\pm 2q^\frac{1}{2}-1$ &
       $q^{19}\pm q^\frac{19}{2}+q^9-q^5+q^2-q\mp q^\frac{1}{2}-1$\\
$22$ & $q^{20}-q^{10}$ & $q^{20}-q^{10}$\\
$23$ & $q^{21}-q^{10}$ & $q^{21}-q^{10}$\\
$24$ & $q^{22}+2q^{11}-q^6\mp 2q^\frac{7}{2}$ & 
       $q^{22}+2q^{11}-q^6\pm q^\frac{7}{2}$\\
$25$ & $q^{23}-q^{11}-q^3+q$ & $q^{23}-q^{11}-q^3+q$\\
$26$ & $q^{24}-q^{12}$ & $q^{24}-q^{12}$\\
$27$ & $q^{25} \mp 2 q^\frac{25}{2}+q^{12}-q^7\pm 2 q^\frac{7}{2}-q^3$ & 
       $q^{25} \pm q^\frac{25}{2}+q^{12}-q^7\mp q^\frac{7}{2}-q^3$\\
$28$ & $q^{26}\pm 2 q^\frac{27}{2}-q^2\mp 2 q^\frac{3}{2}$ & 
       $q^{26}\mp q^\frac{27}{2}-q^2\pm q^\frac{3}{2}$\\
$29$ & $q^{27}-q^{13}$ & $q^{27}-q^{13}$\\
$30$ & $q^{28}\mp2 q^\frac{29}{2}+q^{14}-q^8\pm 2 q^\frac{5}{2}-q^2$ & 
       $q^{28}\pm q^\frac{29}{2}+q^{14}-q^8\mp q^\frac{5}{2}-q^2$
\end{tabular}
\end{center}
where $\pm=(-1)^\frac{r}{2}$.
\end{cor}

\begin{rem}
The expressions for $G_{c,d}(m)$ are approximately of the form 
$q^{m-2}+O(q^\frac{m}{2})$. In the case $2\nmid r$ and $c\ne0$ the formulas for
$G_{c,3}(m)$ are quite close to $q^{m-2}$ when $m=4s$ and $s$ is an odd prime. Indeed, if
$m=4s$ and $s=n^j$ for some odd prime $n$, we see by Corollary \ref{c:nperm},
Theorem \ref{t:expd3}, and Lemma \ref{l:n_elem} that $H_{c,3}(m)=q^{m-2}$, and
then, by Theorem~\ref{thm:Gc}, we have $G_{c,3}(m)=q^{m-2}-q^{\frac{m}{n}-2}$. Especially,
if $j=1$ then $G_{c,3}(m)=q^{m-2}-q^2$.
\end{rem}

\section{Acknowledgment}

This work was inspired by the \emph{Polynomials over Finite Fields and
Applications} workshop at Banff International Research Station, Canada. We
thank the organizers for the invitation to the BIRS workshop.


\begin{thebibliography}{20}

\bibitem{Carlitz52}
L.~Carlitz, 
\newblock A theorem of Dickson on irreducible polynomials, 
\newblock Proc. Amer. Math. Soc. 3 (1952) 693-700.

\bibitem{Carlitz79} 
L.~Carlitz, 
\newblock Explicit evaluation of certain exponential sums, 
\newblock Math. Scand. 44 (1979) 5--16.

\bibitem{Chou01}
W.~S.~Chou, S.~D.~Cohen, 
\newblock Primitive elements with zero traces, 
\newblock Finite Fields Appl. 7 (2001) 125-141.

\bibitem{Cohen00}
S.~D.~Cohen, 
\newblock Kloosterman sums and primitive elements in Galois fields, 
\newblock Acta Arith. 94 (2000) 173-201.

\bibitem{Cohen05}
S.~D.~Cohen,
\newblock Explicit theorems on generator polynomials,
\newblock Finite Fields Appl. 11 (2005) 337-357.

\bibitem{Garcia91}
A.~Garcia, H.~Stichtenoth,
\newblock Elementary abelian $p$-extension of algebraic function fields,
\newblock Manuscripta Math. 72 (1991) 67--79.

\bibitem{Geer92}
G.~van der Geer, R.~Schoof, M.~van der Vlugt,  
\newblock Weight formulas for ternary Melas codes,
\newblock Math. Comp. 58 (1992) 781--792.

\bibitem{Lidl93}
R.~Lidl, G.~L.~Mullen, G.~Turnwald, 
\newblock Dickson Polynomials, 
\newblock Vol. 65 of Pitman Monographs and Surveys in Pure and Applied Mathematics, 
Longman Scientific \& Technical, Harlow, 1993.

\bibitem{Lidl97}
R.~Lidl, H.~Niederreiter,
\newblock Finite Fields,
\newblock Cambridge University Press, Cambridge, 1997.

\bibitem{Moisio06b}
M.~Moisio, 
\newblock The moments of a Kloosterman sum and the weight distribution of a Zetterberg type binary cyclic code, 
\newblock IEEE Trans. Inform. Theory, to Appear. Available: \textsf{http://www.uwasa.fi/\textasciitilde mamo}

\bibitem{Moisio06c}
M.~Moisio, 
\newblock On the moments of Kloosterman sums and fibre products of Kloosterman curves, Submitted. Available: \textsf{http://www.uwasa.fi/\textasciitilde mamo}

\bibitem{Niederreiter90}
H.~Niederreiter,
\newblock An enumeration formula for certain irreducible polynomials
with an application to the construction of irreducible polynomials over the
binary field,
\newblock Appl. Algebra Engrg. Comm. Comput. 1 (1990) 119-124.

\bibitem{Schoof95}
R.~Schoof,
\newblock Families of curves and weight distributions of codes,
\newblock Bull. Am. Math. Soc. 32 (1995) 171--183.

\bibitem{Schoof91}
R.~Schoof, M.~van der Vlugt,
\newblock Hecke operators and the weight distributions of certain codes,
\newblock J. Combin. Theory Ser. A 57 (1991) 163--186.

\bibitem{Stichtenoth} 
H.~Stichtenoth, 
\newblock Algebraic Function Fields and Codes,
\newblock Springer, Berlin, 1993.

\bibitem{Yucas06}
J.~L.~Yucas, 
\newblock Irreducible polynomials over finite fields with prescribed trace/prescribed constant term,
\newblock Finite Fields Appl., in Press.

\end{thebibliography}
\end{document}